\documentclass[preprint]{imsart}

\usepackage{amsfonts}
\usepackage{amssymb}

\usepackage{graphicx}
\usepackage{color}
\usepackage{hyperref}
\usepackage{enumitem}
\usepackage{float}

\RequirePackage[OT1]{fontenc}
\RequirePackage{amsthm,amsmath}
\RequirePackage[numbers]{natbib}
\RequirePackage[colorlinks,citecolor=blue,urlcolor=blue]{hyperref}

\newcommand{\tr}{\operatorname{tr}}
\newcommand{\rank}{\operatorname{rank}}
\newcommand{\diag}{\operatorname{diag}}

\theoremstyle{plain}

\newtheorem{theorem}{Theorem}[section]

\startlocaldefs
\numberwithin{equation}{section}
\theoremstyle{plain}

\endlocaldefs

\begin{document}

\begin{frontmatter}
\title{Reforming the Wishart characteristic function}
\runtitle{Reforming the Wishart characteristic function}

\begin{aug}
\author{\fnms{Eberhard} \snm{Mayerhofer}},

\runauthor{E. Mayerhofer}

\affiliation{University of Limerick}

\address{University of Limerick, Department of Mathematics and Statistics, \\Castletroy, County Limerick, T94 T9PX Ireland. \ead[label=e1]{eberhard.mayerhofer@ul.ie}}
\printead{e1}
\end{aug}

\begin{abstract}
The literature presents the characteristic function of the Wishart distribution on $m\times m$ matrices as an inverse power of the determinant of the Fourier variable, the exponent being
the positive, real shape parameter.

I demonstrate that only for $2\times 2$ matrices, this expression is unambiguous -- in this case the complex range of the determinant excludes the negative real line. When $m\geq 3$, the range of the determinant contains closed lines around the origin, hence a single branch of the complex logarithm does not suffice to define the determinant's power. To resolve this issue, I give the correct analytic extension of the Laplace transform, by exploiting the Fourier-Laplace transform of a Wishart process.
\end{abstract}



\end{frontmatter}

\section{Introduction and Statement of Theorem}
One of the best known multivariate distributions is the Wishart distribution $p(d\xi;\alpha)$, explored first by Wishart in 1928 as
the distribution of the covariance matrix of a normally distributed sample \cite{Wishart}:
\begin{enumerate}
\item For $\alpha$ being a positive half-integer $n/2$, it is the distribution of the covariance of a sample of size $n$ from the $m$-variate
standard normal distribution.\footnote{Using standard normal distributions (zero mean, zero correlations) amounts to the standardization with a scale (matrix) parameter $V=I_m$ 
and avoids a series expansions of the density in terms of zonal polynomials.} More explicitly, let $x_1,\dots, x_n$ be an i.i.d. sequence of
standard normally distributed $m$ vectors. Then the sum of (almost surely) rank one matrices
\begin{equation}\label{eq: qconstr}
x_1x_1^\top+\dots x_{n} x_{n}^\top,
\end{equation}
is Wishart distributed with shape parameter $n/2$. 
\item \label{b} For shape parameters $\alpha>(m-1)/2$, it can be defined in terms of its density\footnote{All notation is summarized in the end of the section.},
\[
p(\xi;\alpha)=\frac{\det(\xi)^{\alpha-\frac{m+1}{2}}e^{-\text{tr}(\xi/2)}}{2^{m\alpha}\Gamma_m(\alpha)}.
\]
\end{enumerate}

The above range of shape parameters $\alpha$ is indeed maximal (see, e.g., \cite{peddada1991proof, EM3, gmm, letac}); the Wishart distribution exists, if 
and only if $\alpha$ lies in
\[
\Lambda_m:=\left\{\frac{1}{2},\dots, \frac{m-2}{2}\right\}\cup \left[\frac{m-1}{2}\right).
\]
This set is typically referred to as Gindikin ensemble (in honor of Gindikin
who discovered this parameterization in the more general setting of homogeneous cones (\cite{Gindikin}).


The characteristic function $\int_{S_m}e^{i\tr(\xi v)}p(d\xi;\alpha)$ of the Wishart distribution is typically\footnote{In his introductory paper \cite{Wishart}, Wishart himself uses densities rather
than characteristic functions.}
quoted as
the ``map"
\begin{equation}\label{eq: det}
\Phi_\alpha: S_m\rightarrow\mathbb C: v\mapsto\det(I_m-2iv)^{-\alpha}.
\end{equation}
Many (standard) references in multivariate statistics state and prove this formula, at least for positive half-integers $\alpha$, e.g. Anderson (18989 citations)\footnote{Google Scholar Data.}, who claims that {\it it can be shown that formula \eqref{eq: det} holds, in general} (that is, for appropriate Fourier-Laplace variables, see \cite[7.3.1, eq. (10) and the subsequent paragraph]{Anderson}), Muirhead (5127 citations) \cite[Theorem 10.3.3]{muirhead2009aspects}, Gupta (1306 citations) \cite[Theorem 3.3.7]{Gupta}, Eaton (839 citations) \cite[Proposition 8.3 (iii)]{eaton} and also the already mentioned work \cite{peddada1991proof}.

To my knowledge, the literature does not offer a rigorous definition of the determinant's power in \eqref{eq: det} (e.g. through the complex power). Most of the references use explicitly
the product formula of the characteristic function to reduce the computation of the characteristic function to a sample of size one: To obtain the formula for
general half-integer $\alpha=n/2$, they thus compute
\[
\left(\det(I_m-2iv)^{-\frac{1}{2}}\right)^n=\det(I_m-2iv)^{-\alpha}
\]
thereby overlooking that the exponential has a complex period, and therefore
the general complex power does in general {\it not} satisfy the functional equation $z^\alpha z^\beta=z^{\alpha+\beta}$. This brief note exemplifies the issues arising
with the formula \eqref{eq: det}, and it provides a clean interpretation for $\Phi_\alpha$, without using the complex power:
\begin{theorem}\label{th main}
The characteristic function of the Wishart distribution with shape parameter $\alpha\in \Lambda_m$
is given by
\[
\int_{S_m^+}e^{i\tr(v\xi)}p(d\xi;\alpha)=e^{\alpha\int_0^1 \text{tr}((I_m-2t iv)^{-1}(2iv))dt},\quad v\in S_m.
\]
\end{theorem}
In the following, I summarize some notation used in this paper, and compute explicitly some characteristic functions, to demonstrate the issue. These are followed
by a proof of the Theorem using the Fourier-Laplace transform of Wishart processes, and some final remarks.
\subsection{Notation}
In this paper, $\mathbb R$ denotes the real line, $i=\sqrt{-1}$ the imaginary unit. For integers $m\geq 1$
I use the following matrix spaces: $S_m$, the real valued symmetric $m\times m$ matrices, $S_m^+$ the positive semidefinite ones,
and $S_m^{++}$ the positive definite matrices. $I_m$ is the unit $m\times m$ matrix, and for a given matrix $A$, $\det(A)$ denotes its determinant and $\tr(A)$ its trace. 

The main branch of the logarithm of a complex number $z$ is defined implicitly by $z=r e^{i\varphi}$, where $\varphi\in (-\pi,\pi]$), in other words,
$\log(z)=\log(r)+i\varphi$.

\section{The problem emerges with size}
No issues arise for $m=1$, where the Wishart distribution coincides with the Gamma distribution with characteristic function
\[
(1-iv\theta)^{-k},
\] 
where $\theta=2$ and $k=\alpha$. In this case it is obvious that $1-2iv$ has only strictly positive real part, whence for general real shape parameter $k$, the
characteristic function should be understood by using the main branch of the logarithm,
\[
(1-2iv)^{-k}=e^{-k\log(1-2iv)}.
\]
Similarly, for dimension $m=2$, an explicit computation shows that the range of $\det(I_2-2iv)$ does not contain the ray $(-\infty,1)$\footnote{In fact, it does not contain the open parabola $y^2<4(1-x)$, which can be
shown by computing the determinant using the eigenvalues of $v$.}. Therefore, again, with the main branch of the complex logarithm, 
\[
\det(I_2-2iv)^{-\alpha}:=e^{-\alpha \log(\det(I_2-2iv))}.
\]

But for any dimension $m\geq 3$, the range of the map $\det(I_m-2iv)$ contains the curve
\[
c: [-\sqrt{3},\sqrt{3}],\quad c(t):=(1-it)^3
\]
(set $2v=\diag(t,t,t,0,\dots,0)$). This curve starts in the complex plane at $(-8,0)$, goes through third and fourth quadrant, intersecting the real axis at $(1,0)$, and continues as its mirror image through first and second quadrant, terminating at its starting point $(-8,0)$.

Therefore, the definition 
\begin{equation}\label{eq: super}
\det(I_m-iv)^{-\alpha}:=\exp(-\alpha \log\det(I_m-2iv)),\quad v\in S_m
\end{equation}
of the characteristic function is ambiguous and a-priori it is not clear that it is complex analytic, because the power is only analytic on a simply connected domain. 

Making the above example more concrete, let me use the main branch of the complex logarithm to define $z^\alpha=e^{\alpha\log(z)}$.
Then I get for $t=\pm \sqrt{3}$, $(1\pm \sqrt{3}i)^3=-8$. Therefore, for $v_\mp=\mp \frac{\sqrt{3}}{2}I_3$
and $\alpha=1/2$,
\begin{equation}\label{eq: 4power}
\det(I_3-2iv_\mp)^{1/2}=((1\pm i\sqrt{3})^3)^{1/2}=\sqrt{-8}=2\sqrt{2} (e^{i\pi})^{1/2}=2\sqrt{2} i.
\end{equation}
But this should be the reciprocal value of the characteristic function of the  sum of three independent gamma distributed random variables 
(the diagonal elements of a Wishart distributed random matrix), hence be of the form
\[
((1\pm i\sqrt{3})^{1/2})^3= (\sqrt 2e^{\pm i\pi/3})^3=\pm 2\sqrt{2} i
\]
which obviously does only agree with \eqref{eq: 4power} for the Fourier variable $v_-$. So for different Fourier variable one needs different branches to get the correct value of
the characteristic function. Similar examples can be produced for arbitrary dimension $m\geq 4$ and for other values of $\alpha$.

Of course, the right side of \eqref{eq: super} is well-defined, if one is willing to restrict the values of the Fourier variable to be of sufficiently small size, relative to $\alpha$, because in a complex neighborhood of the identity matrix, the power of the determinant will be complex analytic.

Beyond that, as suggested by the above examples, the expression $\det(I_m-iv)^{-\alpha}$ could be interpreted by using the eigenvalues of the Fourier variable $u$, due to the fact that the determinant is invariant under actions of the orthogonal group. But this approach will lead to a fairly messy definition. The literature has done computation along these lines (see \cite{Nuttal} and the references therein). The above examples challenge the results of \cite{Nuttal}, however. Besides, the complex analyticity of the characteristic function - as is proved below - cannot easily be read off from a formula that uses a multitude of branches of the complex logarithm.

The next section resolves the problem, by proving the semi-explicit formula for the characteristic function stated in Theorem \ref{th main}.

\section{Proof of Theorem \ref{th main}}\label{sec resolving}
There is an elegant way of interpreting the determinant formula which is inspired by the theory of Wishart processes (introduced by \cite{Bru}, extended to a more general context in \cite{CFMT, CKMT}).
The Wishart process $X=(X_t)_{t\geq 0}$ (here without linear drift term) is the solution to the stochastic differential equation (SDE)
\[
dX_t=\sqrt{X_t}dB_t+dB_t^\top \sqrt{X_t}+2\alpha I_m dt,\quad X_0=x\in S_m^+,
\]
where all products are in terms of the matrix multiplication, and $B$ is an $m\times m$ matrix of standard Brownian motions. This SDE
indeed has a unique global weak solution if $\alpha\geq (m-1)/2$, or ($\rank(x)\leq 2\alpha$, and $\alpha\in \{1/2,\dots,(m-2)/2\}$), see \cite[Theorem 1.3]{gmm}. Furthermore, by
\cite[Proposition 2.8 and Lemma 2.9]{EM2} for any $u\in S_m^+$
\[
\mathbb E[e^{-\tr(uX_t)}]=e^{-\phi(t,u)-\tr(\psi(t,u)x)},\quad t\geq 0,
\]
where the matrix valued function $\psi(t,u)$ satisfies the initial value problem
\begin{align*}
\partial_t\psi(t,u)=-2\psi(t,u)^2,\quad \psi(t=0,u)=u,
\end{align*}
and thus is of the explicit form $\psi(t,u)=(I_m+2t u)^{-1}u$. 
The function $\phi(t,u)$ satisfies the trivial equation 
\[
\partial_t\phi(t,u)=2\alpha\text{tr}(\psi(t,u)),\quad \phi(t=0,u)=0,
\] 
and thus equals 
\begin{equation}\label{eq: phi lt}
\phi=\alpha\log(\det(I_m+2tu)).
\end{equation}
In particular, for $x=0$, the Wishart SDE has unique global weak solution
\eqref{eq: phi lt}, hence it satisfies at $t=1$
\[
e^{-\phi(1,u)}=\det(I_m+2u)^{-\alpha}.
\]
This is exactly the Laplace transform of the Wishart distribution.

Consider now the (trivial) initial value problem for the function
$\phi(t,u-iv): [0,1]\times \mathbb D\rightarrow \mathbb C$, where
\[
\mathbb D=\{u-iv\in S_m+iS_m\mid u\in -I_m/2+S_m^{++}\}.
\]
\begin{equation}\label{eq: diff eq}
\partial_t\phi(t,u-iv)=\alpha\text{tr}\left((I_m+2t(u-iv))^{-1}(2u-2iv)\right),\quad \phi(t=0,u-iv)=0.
\end{equation}
The IVP defined in \eqref{eq: diff eq} has also unique solutions for any $u\in \mathbb D$ with blow-up beyond $t=1$, because $(I_m+2t(u-iv))$ is non-singular for each $t\in [0,1]$.\footnote{By Ostrowski and Tausski (cf. \cite[(4.4)]{Johnson}, a complex matrix $X+iY$, where $X,Y$ are real and $X$ is positive definite, satisfies $\vert\det(X+iY)\vert\geq \det(X)$.} 
Therefore, the solution of \eqref{eq: diff eq} can be written semi-explicitly as
\[
\phi(1,u+iv):=\alpha\int_0^1 \text{tr}((I_m+2t (u-iv))^{-1}(2u-2iv))dt.
\]
This solution is complex analytic in $d(d+1)/2$ variables (the entries of $u-iv$), because the right-side of \eqref{eq: diff eq} is (cf. \cite[Theorem 10.8.2]{Dieudonne}). Hence
\[
e^{-\phi(1,u-iv)}
\]
constitutes an analytic extension of $\det(I_m+2u)^{-\alpha}$ to the complex strip $\mathbb D$. Now 
$e^{-\phi(1,u-iv)}$ and $\int_{S_m^+}e^{-\tr((u-iv)\xi)}p(\xi;\alpha)d\xi$ are both complex analytic in several variables on $\mathbb D$, and the two functions agree for $v=0$ and $u\in S_m^{++}$. An open domain in $S_m$
is a set of uniqueness, hence the two functions agree (cf. \cite[(9.4.4)]{Dieudonne}).

The proof of Theorem \ref{th main} is complete.

\section{Final remarks}

The results of this paper can easily be extended to the non-central Wishart distribution, the difference being a multiplicative exponential term in the characteristic
function (instead of $x=0$ use in Section \ref{sec resolving} a general starting point $X_0=x$ for the Wishart process, respecting the rank condition, see \cite{gmm}).

Analytic continuation arguments similar to the above are extensively used in the recent literature on affine Markov processes, e.g. in the papers
\cite{CFMT, CKMT, DFEM, MKREM}. 


\bibliographystyle{amsplain}

\end{document}